\documentclass[12pt]{amsart}
\usepackage{graphicx}
\usepackage{amsmath,amssymb}
\newtheorem{theorem}{Theorem}[section]

\newtheorem{lemma}[theorem]{Lemma}

\newtheorem{corollary}[theorem]{Corollary}

\theoremstyle{remark}

\numberwithin{equation}{section}

\begin{document}

\title[
Unknotting rectangular diagrams of the trivial knot
]{
Unknotting rectangular diagrams of the trivial knot by exchanging moves
}

\author{Chuichiro Hayashi and Sayaka Yamada}

\date{\today}

\thanks{The first author is partially supported
by JSPS KAKENHI Grant Number 22540101.}

\begin{abstract}
 If a rectangular diagram represents the trivial knot,
then it can be deformed into the rectangular diagram with only two vertical edges
by a finite sequence of merge operations and exchange operations,
without increasing the number of vertical edges,
which was shown by I. A. Dynnikov.
 We show in this paper that we need no merge operations
to deform a rectangular diagram of the trivial knot
to one with no crossings.
\end{abstract}

\maketitle

\section{Introduction}\label{sect:introduction}
 Birman and Menasco introduced arc-presentation of knots and links in \cite{BM},
and Cromwell formulated it in \cite{C}.
 Dynnikov pointed out in \cite{D1} and \cite{D2} 
that Cromwell's argument in \cite{C} almost shows 
that any arc-presentation of a split link
can be deformed into one which is $\lq\lq$visibly split"
by a finite sequence of exchange moves. 
 He also showed
that any arc-presentation of the trivial knot
can be deformed into one with only two arcs
by a finite sequence of merge moves and exchange moves, 
without using divide moves which increase the number of arcs.
 As is shown in page 41 in \cite{C},
an arc-presentation is equivalent to a rectangular diagram.
 All arguments in this paper are described in words on rectangular diagrams.

\begin{figure}[htbp]
\begin{center}
\includegraphics[width=50mm]{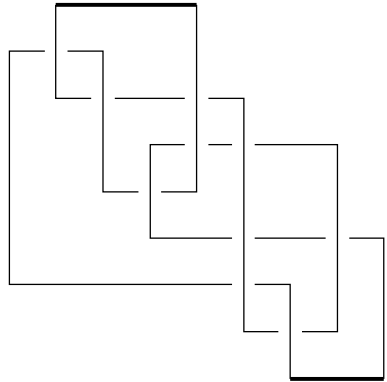}
\end{center}
\caption{rectangular diagram}
\label{fig:TrivialKnot9arcs}
\end{figure}

 A {\it rectangular diagram} of a knot is a knot diagram
in the plane ${\mathbb R}^2$
which is composed of vertical lines and horizontal lines 
such that no pair of vertical lines are colinear,
no pair of horizontal lines are colinear, 
and the vertical line passes over the horizontal line at each crossing. 
 See Figure \ref{fig:TrivialKnot9arcs}.
 These vertical lines and horizontal lines 
are called {\it edges} of the rectangular diagram.
 Every knot or link has a rectangular diagram 
(Proposition in page 42 in \cite{C}).

\begin{figure}[htbp]
\begin{center}
\includegraphics[width=90mm]{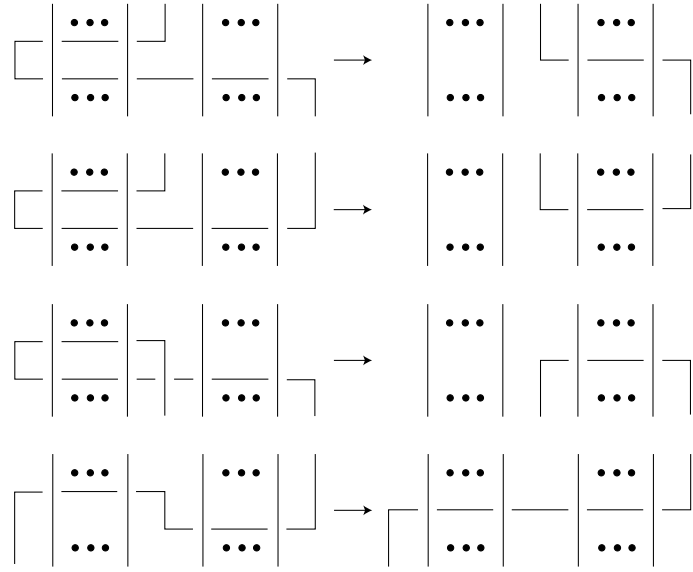}
\end{center}
\caption{interior horizontal merge}
\label{fig:merge1}
\end{figure}

 Cromwell moves,
which are described in the next three paragraphs,
are elementary moves for rectangular diagrams of knots and links.
 They do not change type of knots and links.
 Moreover, Theorem in page 45 in \cite{C} and Proposition 4 in \cite{D1}
state that, if two rectangular diagrams represent the same knot or link,
then one is obtained from the other 
by a finite sequence of these elementary moves
and rotation moves, which is introduced below.

\begin{figure}[htbp]
\begin{center}
\includegraphics[width=70mm]{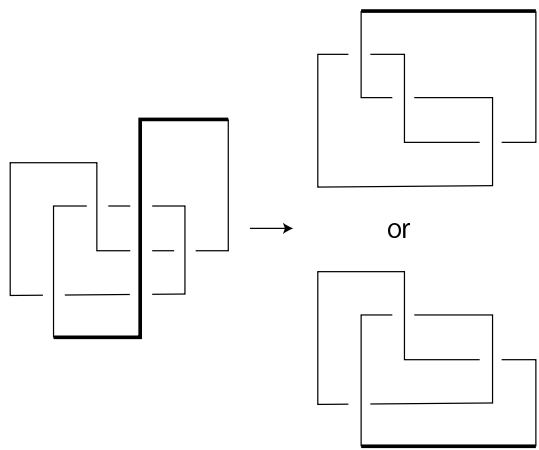}
\end{center}
\caption{exterior horizontal merge}
\label{fig:merge2}
\end{figure}

 First, we recall merge moves.
 If two horizontal (resp. vertical) edges 
connected by a single vertical (resp. horizontal) edge
have no other horizontal (resp. vertical) edges between their abscissae (resp. ordinates),
then we can amalgamate the three edges
into a single horizontal (resp. vertical) edge.
 This move is called an {\it interior horizontal (resp. vertical) merge}.
 See Figure \ref{fig:merge1} for examples of interior horizontal merge moves.
 If two horizontal (resp. vertical) edges 
connected by a single vertical (resp. horizontal) edge
have the other horizontal (resp. vertical) edges between their abscissae (resp. ordinates),
i.e., they are the top and bottom (the leftmost and rightmost) edges,
then we can amalgamate the three edges
into a single horizontal (resp. vertical) edge.
 We may place the new horizontal (resp. vertical) edge
either at the top height or at the bottom height
(resp. either in the leftmost position or in the rightmost position).
 See Figure \ref{fig:merge2}.
 We call this move an {\it exterior horizontal (resp. vertical) merge}.
 Note that a merge move decreases 
the number of vertical edges and that of horizontal edges by one.
 The inverse moves of merge moves are called {\it divide moves}.

\begin{figure}[htbp]
\begin{center}
\includegraphics[width=90mm]{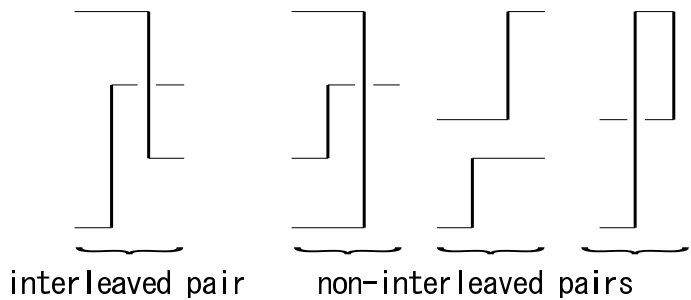}
\end{center}
\caption{interleaved pair and non-interleaved pairs}
\label{fig:interleaved}
\end{figure}

 To describe exchange moves, we need a terminology.
 Two vertical edges are said to be {\it interleaved},
if the heights of their endpoints alternate.
 See Figure \ref{fig:interleaved}.
 Similarly, we define interleaved two horizontal edges.

\begin{figure}[htbp]
\begin{center}
\includegraphics[width=90mm]{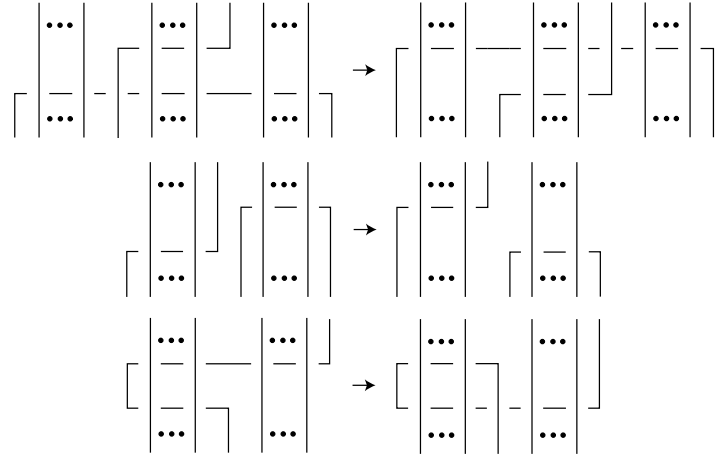}
\end{center}
\caption{interior horizontal exchange}
\label{fig:exchange1}
\end{figure}

\begin{figure}[htbp]
\begin{center}
\includegraphics[width=50mm]{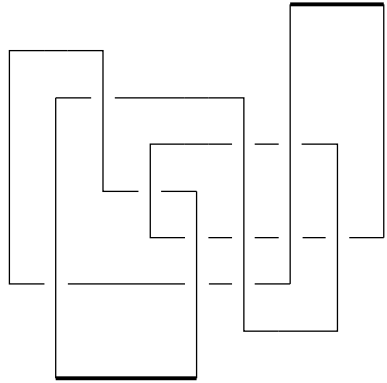}
\end{center}
\caption{This is obtained from Figure \ref{fig:TrivialKnot9arcs} by an exterior horizontal exchange.}
\label{fig:exchange2}
\end{figure}

 If two horizontal edges at mutually adjacent heights are not interleaved,
then we can exchange their heights.
 See Figure \ref{fig:exchange1}.
 This move is called an {\it interior horizontal exchange}.
 If the top horizontal edge and the bottom one are not interleaved,
then we can exchange their heights.
 We call this move an {\it exterior horizontal exchange}.
 See Figure \ref{fig:exchange2},
where the rectangular diagram obtained from one in Figure \ref{fig:TrivialKnot9arcs}
by an exterior horizontal exchange move.
 Similarly, we define {\it vertical exchange} moves.

 A {\it horizontal rotation move} on a rectangular diagram
moves the top edge to the bottom, or the bottom edge to the top.
 A {\it vertical rotation move} on a rectangular diagram 
moves the leftmost edge to the right, or the rightmost edges to the left.
 (The horizontal (resp. vertical) rotation move corresponds 
to the $\pm 2\pi/n$ rotation about the dual axis (resp. the axis)
on an arc-presentation.)
 However, we do not use rotation move in this paper.

 The next result of Dynnikov gives a finite algorithm
to decide whether a given rectangular diagram represents
the trivial knot or not.

\begin{theorem}[Dynnikov \cite{D1}, \cite{D2}]
 Any rectangle diagram of the trivial knot
can be deformed into one with only two vertical edges
by a finite sequence of merge moves and exchange moves.
\end{theorem}

 Note that the sequence in the above theorem contains no divide moves.
 Hence the sequence gives a monotone simplification,
that is, no move in the sequence increases the number of vertical edges.
 There are only finitely many rectangle diagrams
with a fixed number of vertical edges.
 Thus the above theorem gives a finite algorithm
for the decision problem.

 The sequence as in Dynnikov's theorem 
sometimes needs to contain exterior exchange moves. 
 In fact, the rectangle diagram shown in Figure \ref{fig:TrivialKnot9arcs}
represents the trivial knot.
 It admits no merge moves
since it does not have an edge of length $1$ or $9-1$.
 We cannot apply any interior horizontal exchange move to the diagram
because every pair of horizontal edges in adjacent levels
are interleaved.
 It can be seen easily
that no vetical exchange move can be performed on this diagram.
 Hence every sequence as in Dynnikov's theorem on this diagram
must begin with the exetrior horizontal exchange move.

 In \cite{HK}, A. Henrich and L. Kauffman announced
an upper bound of the number of Reidemeister moves needed for unknotting
by applying Dynnikov's theorem to rectangular diagrams. 
 Lemma 7 in \cite{HK} states 
that no more than $n-2$ Reidemeister moves are required
to perform an exchange move on a rectangular diagram with $n$ vertical edges.
 However, their proof of Lemma 7 in \cite{HK} does not consider
the exterior exchange moves.

 In this paper,
we show that the decision problem can be solved 
without using merge moves.
 In fact, we obtain the next result 
by taking advantage of Dynnikov's theorem.

\begin{theorem}\label{theorem:main}
 For any pair of rectangle diagrams of the trivial knot
with the same number of vertical edges,
there is a finite sequence of exchange moves
which deform one into the other.
\end{theorem}

\begin{corollary}
 Any rectangle diagram of the trivial knot
can be deformed into one with no crossings
by a finite sequence of exchange moves.
\end{corollary}

 We need much larger number of Cromwell moves 
in a sequence as in this corollary
than in that in Dynnikov's theorem.
 However, this improves Lemma 6 in \cite{HK}
which assures that
the number of exchange moves and merge moves
in a sequence as in Dynnikov's theorem
is bounded above by $\sum_{i=2}^n \frac{1}{2}i[(i-1)!]^2$,
where $n$ is the number of vertical edges
of the first rectangular diagram.
 In fact, we obtain the upper bound $\frac{1}{2}n[(n-1)!]^2$
which is an upper bound 
for the number of combinatorially distinct rectaugluar diagrams with $n$ vertical edges
given in Proposition 5 in \cite{HK}

\section{Proof of Theorem \ref{theorem:main}}\label{section:proof}

 We prove Theorem \ref{theorem:main} in this section.

\begin{figure}[htbp]
\begin{center}
\includegraphics[width=100mm]{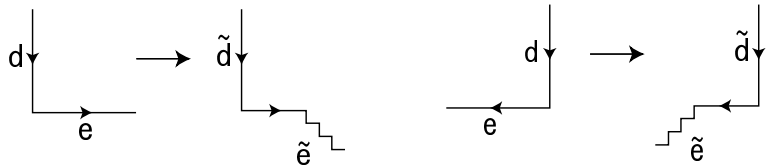}
\end{center}
\caption{nothced edge}
\label{fig:NotchedEdge}
\end{figure}

\begin{figure}[htbp]
\begin{center}
\includegraphics[width=70mm]{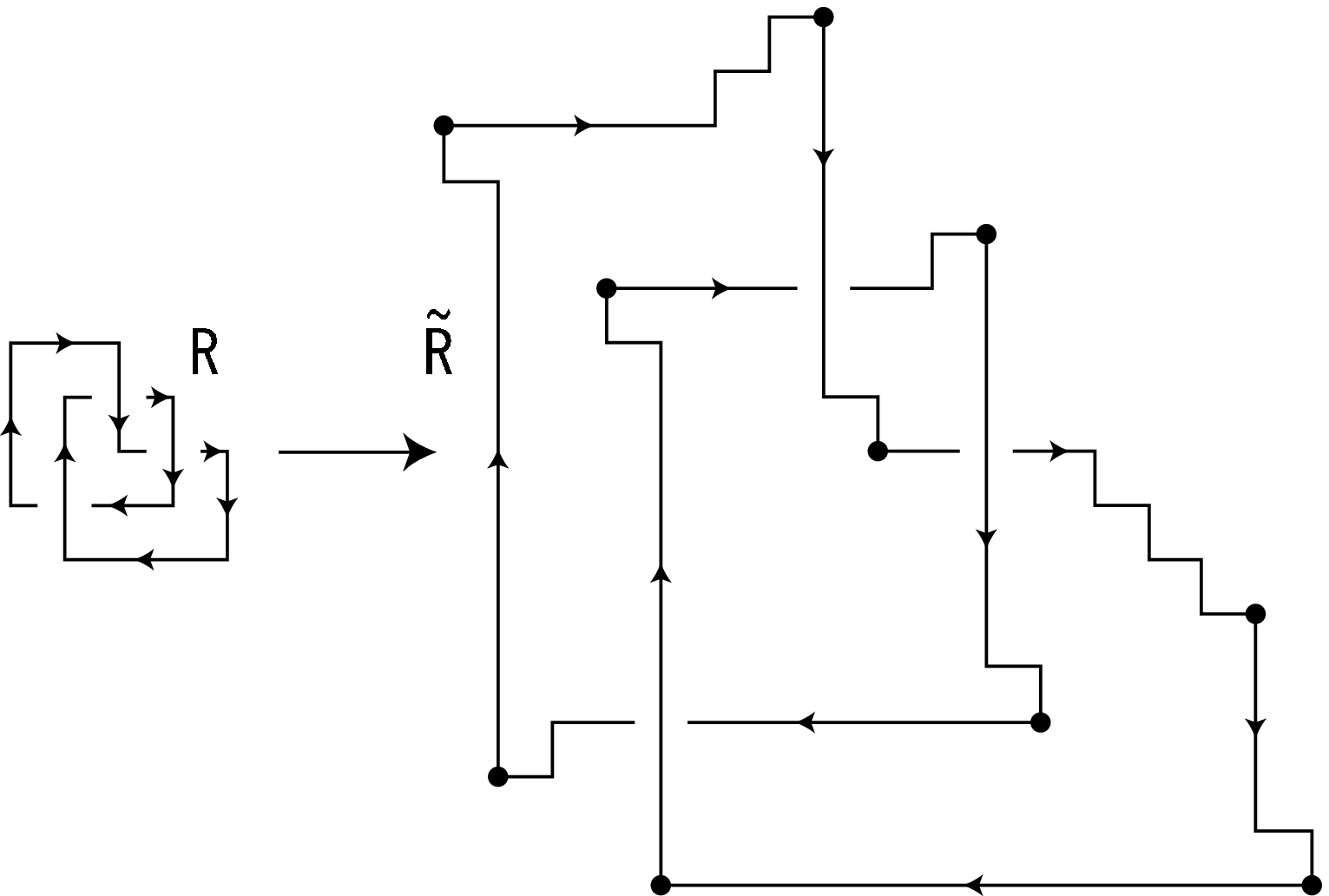}
\end{center}
\caption{outline}
\label{fig:outline}
\end{figure}

 Let $R$, $\tilde{R}$ be oriented rectangular diagrams of the same knot.
 We say that $R$ is an {\it outline} of $\tilde{R}$
if $\tilde{R}$ is obtained from $R$
by replacing the edges of $R$
by {\it notched edges} as below.
 Let $e$ be an original edge of $R$,
and $d$ the edge preceding $e$.
 Then the notched edge $\tilde{e}$ substituting for $e$
consists of a long edge 
followed by $0$ or even number of very short edges
going away from the right angle between $d$ and $e$
as in Figure \ref{fig:NotchedEdge}.
 We say that a notched edge has $r$ notches if it consists of $2r+1$ edges.
 In Figure \ref{fig:NotchedEdge}, $\tilde{d}$ has no notches for simplicity.
 See Figure \ref{fig:outline} for entire view of $R$ and $\tilde{R}$.
 Note that $\tilde{R}$ has the same number of crossings as $R$.

\begin{lemma}\label{lemma:SingleNotchedEdge}
 Let $R$ be an oriented rectangular diagram of a knot,
and $\tilde{R}$ an oriented rectangular diagram
which has $R$ as its outline.
 Then for each edge $e$ of $R$,
there is an oriented rectangular diagram $\tilde{R}_e$ 
such that (1) $\tilde{R}_e$ is obtained from $\tilde{R}$
by a finite sequence of exchange moves,
that (2) $R$ is an outline of $\tilde{R}_e$
and that (3) $\tilde{R}$ has only one notched edge $\tilde{e}$
with positive number of notches
and $\tilde{e}$ substitutes for $e$.
\end{lemma}

\begin{figure}[htbp]
\begin{center}
\includegraphics[width=80mm]{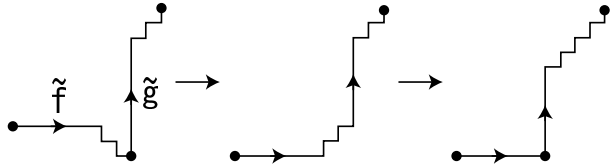}
\end{center}
\caption{moving notches}
\label{fig:MoveNotches}
\end{figure}

\begin{proof}
 Let $S$ be an oriented rectangular diagram of a knot,
and $\tilde{S}$ an oriented rectangular diagram 
which has $S$ as its outline.
 Let $f$ be an edge of $R$,
$g$ the edge of $R$ next to $f$,
$\tilde{f}$ and $\tilde{g}$ notched edges of $\tilde{S}$ corresponding to $f$ and $g$, 
$r_f$ and $r_g$ 
the number of notches of $\tilde{f}$ and $\tilde{g}$ respectively. 
 We show that 
an adequate finite sequence of exchange moves
brings the $r_f$ notches of $\tilde{f}$ 
to behind the tail of the long edge of $\tilde{g}$.

 Assume, without loss of generality,
that $f$ is horizontal and $g$ is vertical.
 Let $c_g$ be the number of horizontal edges which cross $g$.
 If the notches of $\tilde{f}$ goes into the right angle between $f$ and $g$,
then we simply 
exchange horizontal edges, say $h_f$, of notches of $f$
and $c_g$ horizontal notched edges crossing $\tilde{g}$. 
 Precisely, we need more exchange moves
if there are horizontal edges away from $\tilde{g}$ 
between the height of the top of $h_f$ 
and the height of the bottom of notches of $\tilde{g}$.
 If the notches of $\tilde{f}$ 
goes away from the right angle between $f$ and $g$,
then we first exchange every pair of $r_f + 1$ horizontal edges of $\tilde{f}$
to obtain the situation as in the previous case.
 See Figure \ref{fig:MoveNotches}.
 Note that edges of notches are not interleaved with any edge
since they are very short.

 Then, by applying operations as above to $R$ repeatedly,
we obtain the lemma.
\end{proof}

\begin{figure}[htbp]
\begin{center}
\includegraphics[width=60mm]{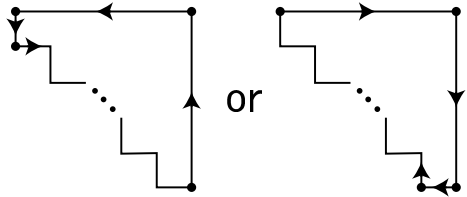}
\end{center}
\caption{the rectangular diagram $\Delta$ with $\Phi$ being an outline}
\label{fig:goal}
\end{figure}

 Let $R_0$ be a rectangular diagram with $n$ vertical edges
which represents the trivial knot.
 We give $R_0$ an arbitrary orientation.
 By Dynnikov's theorem, 
there is a finite sequence of merge moves and exchange moves
which deforms $R_0$ to the rectangle diagram 
$\Phi$ 
with only two vertical edges,
and performs a merge move 
whenever there are applicable merge moves.
 We take such a sequence, fix it, and call it a {\it Dynnikov sequence}.
 Let 
$R_0 \rightarrow R_1 \rightarrow R_2 \rightarrow \cdots \rightarrow R_m=\Phi$
be the sequence of the rectangular diagrams
such that $R_i$ is obtained from $R_{i-1}$
by the $i$\,th move in the Dynnikov sequence.
 We set $\tilde{R}_0 = R_0$,
and will form a rectangle diagram $\tilde{R}_i$ with $R_i$ being an outline
so that $\tilde{R}_{i+1}$ is obtained from $\tilde{R}_i$
by a finite sequence of exchange moves.
 Note that $\tilde{R}_m$ has no crossings
because $R_m$ has only two vertical edges.
 Moreover,
$\tilde{R}_m$ can be deformed
into one of the two diagrams in Figure \ref{fig:goal}
by Lemma \ref{lemma:SingleNotchedEdge}.
 Note that, if we ignore their orientations, 
these two rectangular diagrams are the same,
which we denote by $\Delta$.
 This implies Theorem \ref{theorem:main}.
 Let $T_1$ and $T_2$ be rectangular diagrams representing the trivial knot.
 For $i=1$ and $2$, 
there is a sequence $\gamma_i$ of merge and exchange moves 
which deforms $T_i$ into $\Delta$.
 Then the sequence $\gamma_1$ followed 
by the inverse of $\gamma_2$ deforms $T_1$ into $T_2$.
 Hence, it is enough to show the next lemma.

\begin{lemma}\label{lemma:exchange}
 Let $S_0 \rightarrow S_1 \rightarrow S_2 \rightarrow \cdots \rightarrow S_m$
be a sequence of rectangular diagrams of a knot
such that $S_i$ is obtained from $S_{i-1}$
by a merge or exchange move.
 Then there is another sequence 
$\tilde{S}_0 \rightarrow \tilde{S}_1 \rightarrow \tilde{S}_2 
\rightarrow \cdots \rightarrow \tilde{S}_m$ with $\tilde{S}_0=S_0$
such that $S_i$ is an outline of $\tilde{S}_i$
and that $\tilde{S}_i$ is obtained from $\tilde{S}_{i-1}$ 
by a finite sequence of exchange moves.
\end{lemma}

 In the above lemma,
note that the knot may not be trivial.
 We show the way how to construct $\tilde{S}_i$ from $\tilde{S}_{i-1}$
refering to the $i$\,th move $S_{i-1} \rightarrow S_i$ 
in the original sequence.

%

 First, we apply Lemma \ref{lemma:SingleNotchedEdge} to $\tilde{S}_{i-1}$
so that the notches are gathered to a single notched edge
which is away from the notched edges 
corresponding to the edges participating in the $i$\,th move on $S_{i-1}$. 

 When the $i$\,th move is an exchange move, 
we simply exchange the corresponding notched edges with no notches.

\begin{figure}[htbp]
\begin{center}
\includegraphics[width=35mm]{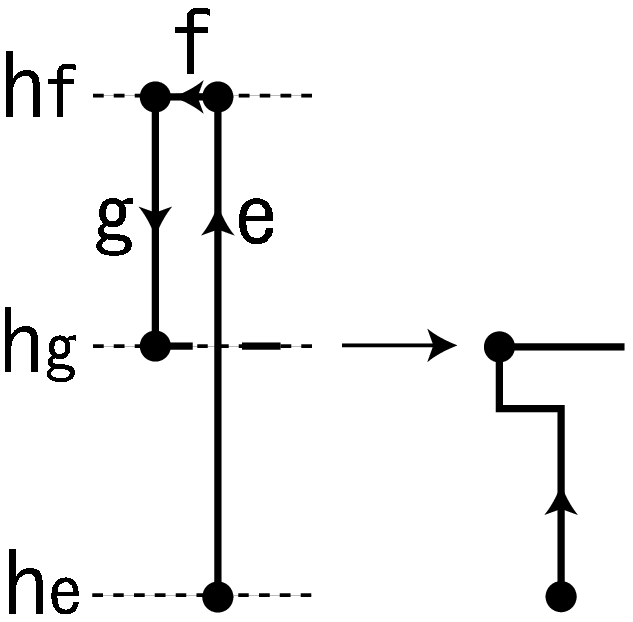}
\end{center}
\caption{a sequence of exchange moves substituting for an interior vertcial merge move}
\label{fig:ExchangeForMerge1}
\end{figure}

 In the case where the $i$\,th move is a merge move,
let $e$, $f$ and $g$ be edges of $S_{i-1}$
such that they appear in succession in this order in the knot
and that the $i$\,th move merges these three edges into a single edge.
 Assume, without loss of generality,
that $e$ and $g$ are vertical and $f$ is horizontal.
 Let $h_f$ be the height of $f$,
and $h_e$ and $h_g$ the heights of endpoints of $e$ and $g$
which are not equal to $h_f$.
 The notched edges of $\tilde{S}_{i-1}$ corresponding to $e,f,g$ have no notches 
and are called $e,f,g$ for simplicity in the followings.

 First, suppose that the $i$\,th move is an interior merge move.
 Then we move $f$ in $\tilde{S}_{i-1}$
to a height $h$ adjacent to $h_g$ by exchange moves.
 When $h_g < h_e$ (resp. $h_e < h_g$),
we take $h$ to be $h_g + \epsilon$ (resp. $h_g - \epsilon$),
where $\epsilon$ is a very small positive real number.
 See Figure \ref{fig:ExchangeForMerge1}.

\begin{figure}[htbp]
\begin{center}
\includegraphics[width=100mm]{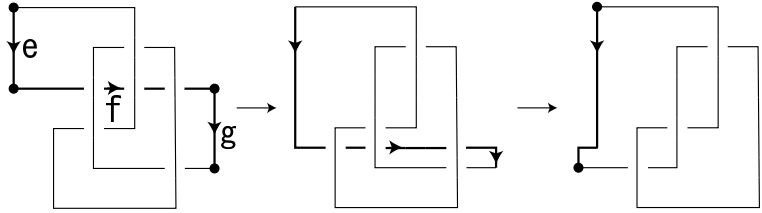}
\end{center}
\caption{a sequence of exchange moves substituting for an exterior vertical merge move}
\label{fig:ExchangeForMerge2}
\end{figure}

\begin{figure}[htbp]
\begin{center}
\includegraphics[width=125mm]{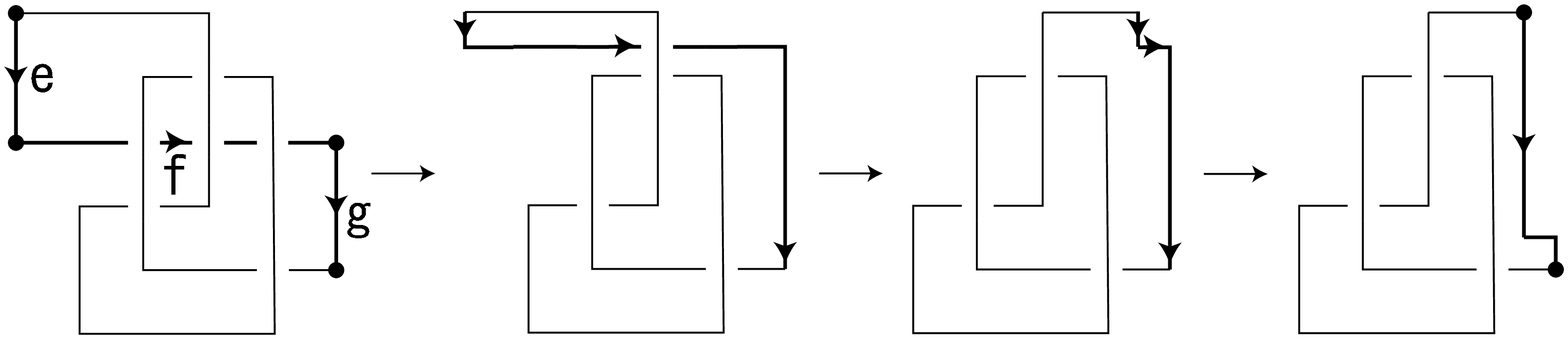}
\end{center}
\caption{a sequence of exchamge moves substituting for an exterior vertical merge move}
\label{fig:ExchangeForMerge3}
\end{figure}

 When the $i$\,th move is an exterior merge move,
let $k$ be the edge of $S_i$ obtained by amalgamating the edges $e, f, g$.
 First, we consider the case 
where the new edge $k$ is in the vertical line of the same lateral coordinate as $e$.
 We move $f$ in $\tilde{S}_{i-1}$
by exchange moves to the same height as in the previous paragraph,
to make $g$ very short.
 Then we move $g$ by exchange moves to the position very close to $e$
so that it forms an edge of a notch.
 In the case where $k$ is in the vertical line of the same lateral coordinate as $g$, 
we move $f$ in $\tilde{S}_{i-1}$
by exchange moves to the height $h_e + \epsilon$ (resp. $h_e - \epsilon$)
to make $e$ very short
when $h_e < h_g$ (resp. $h_e > h_g$).
 Then we move shortened $e$ by exchange move 
to the position slightly before $g$.
 Then we move the horizontal edge between $e$ and $g$ by exchange moves
to the height very close to $h_g$
so that it forms an edge of a notch.
 See Figure \ref{fig:ExchangeForMerge3}.
 This completes the proof of Lemma \ref{lemma:exchange},
and hence of Theorem \ref{theorem:main}.

\bibliographystyle{amsplain}

\medskip

\noindent
Chuichiro Hayashi: 
Department of Mathematical and Physical Sciences,
Faculty of Science, Japan Women's University,
2-8-1 Mejirodai, Bunkyo-ku, Tokyo, 112-8681, Japan.
hayashic@fc.jwu.ac.jp

\vspace{3mm}
\noindent
Sayaka Yamada:
Department of Mathematical and Physical Sciences,
Faculty of Science, Japan Women's University,
2-8-1 Mejirodai, Bunkyo-ku, Tokyo, 112-8681, Japan.
\newline
m1136008ys@gr.jwu.ac.jp

\end{document}